\documentclass{article}

     \PassOptionsToPackage{numbers, compress}{natbib}

     \usepackage[final]{neurips_2023_ml4ps}

\usepackage[utf8]{inputenc} 
\usepackage[T1]{fontenc}    
\usepackage{hyperref}       
\usepackage{url}            
\usepackage{booktabs}       
\usepackage{amsfonts}       
\usepackage{nicefrac}       
\usepackage{microtype}      
\usepackage{xcolor}         

\usepackage{bm}             
\usepackage{graphicx}       

\title{A Data-Driven, Non-Linear, Parameterized Reduced Order Model of Metal 3D Printing}

\author{%
    Aaron L. Brown \\
    Department of Mechanical Engineering, Stanford University
    \\
    Stanford, CA, 94305
    \\
    \texttt{abrown97@stanford.edu}
    \\
    \And
    Eric B. Chin \\
    Lawrence Livermore National Laboratory \\
    Livermore, CA, 94550 \\
    \texttt{chin23@llnl.gov}
    \And
    Youngsoo Choi \\
    Lawrence Livermore National Laboratory \\
    Livermore, CA, 94550 \\
    \texttt{choi15@llnl.gov}
    \And
    Saad A. Khairallah \\
    Lawrence Livermore National Laboratory \\
    Livermore, CA 94550 \\
    \texttt{khairallah1@llnl.gov}
    \And
    Joseph T. McKeown \\
    Lawrence Livermore National Laboratory \\
    Livermore, CA 94550 \\
    \texttt{mckeown3@llnl.gov}
}

\begin{document}

\maketitle

\begin{abstract}
Directed energy deposition (DED) is a promising metal additive manufacturing technology capable of 3D printing metal parts with complex geometries at lower cost compared to traditional manufacturing. The technology is most effective when process parameters like laser scan speed and power are optimized for a particular geometry and alloy. To accelerate optimization, we apply a data-driven, parameterized, non-linear reduced-order model (ROM) called Gaussian Process Latent Space Dynamics Identification (GPLaSDI) to physics-based DED simulation data. With an appropriate choice of hyperparameters, GPLaSDI is an effective ROM for this application, with a worst-case error of about $8\%$ and a speed-up of about $1,000,000$x with respect to the corresponding physics-based data.
\end{abstract}

\section{Introduction}
In Directed Energy Deposition (DED) additive manufacturing, or metal 3D printing, metal powder is fed through a nozzle and directed down to a build plate, where a laser fuses the incoming powder into the previous layer or substrate. 
An important goal in DED is to optimize process parameters such as laser power and scan speed to control the temperature history and prevent defects in prints. 
Currently, the additive manufacturing community resorts to a trial and error approach involving numerous experiments, each of which can be time-consuming and costly.
This process must be repeated for any new printing material, and the results may vary depending on the complexity of the part or the 3D printer itself. 

Previous approaches have attempted analytic solutions to this problem, which in its simplest form asks the temperature field for a moving, distributed heat source \cite{eagar1983temperature, honarmandi2021rigorous, wolfer2019fast}. However, analytic approaches do not offer the flexibility to handle irregular domains and complicated initial and boundary conditions found in many additive manufacturing applications. Such flexibility is afforded by high-fidelity finite element-based models of DED \cite{khairallah2023high}, but simulations can take weeks to run, making them unsuitable for parameter optimization.

Our approach is to use reduced order models (ROMs), informed by simulations or experiments. ROMs can be orders of magnitude faster with a marginal reduction in accuracy and have been successfully applied to various physics problems \cite{carlberg2018conservative, cheung2023datascarce, lauzon2022s, tsai2023accelerating, fries2022lasdi, kim2021efficient, copeland2022reduced, huhn2023parametric, mcbane2022stress, choi2020sns, cheung2023local, hoang2021domain}. They may fully replace full-order simulations or experiments or be used as surrogate models such as in the Surrogate Management Framework \cite{marsden2004optimal}. In this work, we apply the recently-developed ROM known as Gaussian Process Latent Space Dynamics Identification (GPLaSDI) \cite{bonneville2024gplasdi} to parameterized DED simulation data to predict, for given process parameters, the time-dependent temperature field, which may then be used as an input to a microstructure evolution model. A major advantage of GPLaSDI is its use of an autoencoder to spatially compress data in a nonlinear manner, which, compared to linear compression methods such as the singular value decomposition, can better capture the advection-like behavior introduced by the motion of the DED print head.

\section{Methods}
\subsection{DED simulation data}
In this work, we consider data obtained from DED simulations performed in ALE3D \cite{alewebsite}. The data was parameterized with respect to two important parameters in DED: laser power $P$ and laser scan speed $S$. Specifically, we generated single-track DED simulations for each sample on a 5x5 grid in parameter space -- $P \in [120, 130, 140, 150, 160]$ W, $S \in [0.08, 0.09, 0.10, 0.11, 0.12]$ m/s. Each simulation modeled the deposition of a single track of titanium alloy Ti-6Al-4V (Ti64) at a rate of 5 g/min over a length of 2.5mm. The computational mesh consisted of $N_u = 25,990$ nodes, and each simulation was run between $100$ and $150$ ms, from which we extract $N_t + 1 = 101$ evenly spaced time steps for our dataset. Each simulation took on average 2.5 hours on a single core of an Intel(R) Xeon(R) Platinum 8479 CPU.  From this dataset, we select $N_\mu$ samples for the training set; how the training set is chosen will be described later. The remainder are part of the test set. 

\begin{figure}
  \centering
  \includegraphics[width=\textwidth]{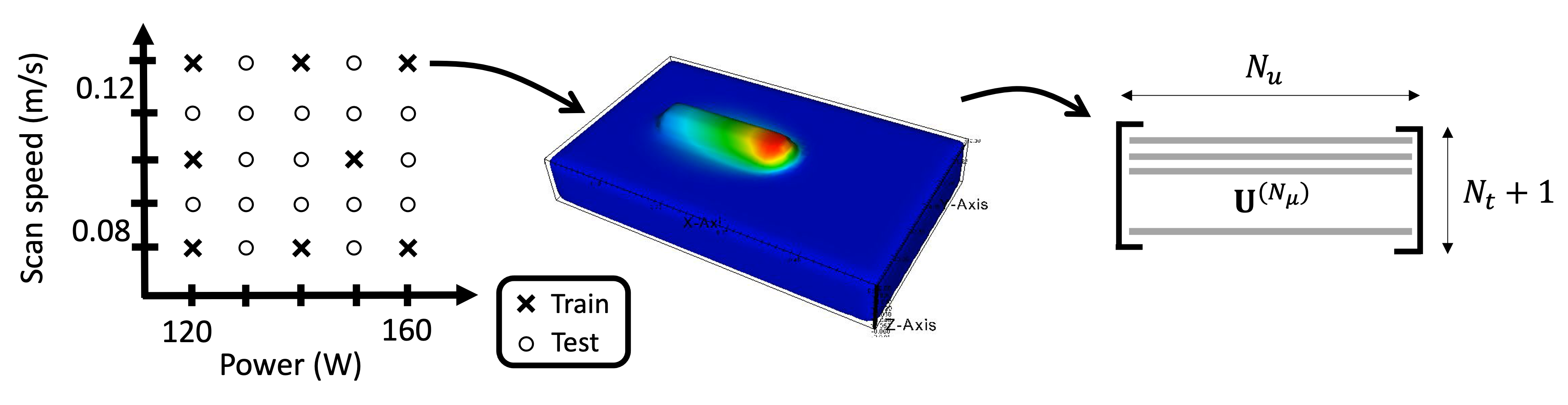}
  \caption{Schematic of dataset generation. Single-track DED simulations were generated for parameters (laser power and laser scan speed) on a 5x5 grid. Temperature data was extracted for each simulation and arranged into a training snapshot matrix.}
  \label{fig:dataset_generation}
\end{figure}

Let $\boldsymbol{\mu}^{(i)} = (P^{(i)}, S^{(i)} )$ be the $i^{th}$ parameter vector in the training set. The quantity of interest from these simulations is the time-dependent temperature field, $u$. Let $\mathbf{u}_n^{(i)} \in \mathbb{R}^{N_u}$ be the snapshot vector of temperature data for all nodes at time step $n$ for parameter vector $\boldsymbol{\mu}^{(i)}$, and let
$
    \mathbf{U}^{(i)} = \Big[ 
    \mathbf{u}_0^{(i)}, ..., \mathbf{u}_{N_t}^{(i)}
    \Big]
    \in
    \mathbb{R}^{(N_t + 1) \times N_u}
$
be the data matrix for all nodes and all time steps for parameter vector $\boldsymbol{\mu}^{(i)}$. Combining all these matrices together, we define a 3rd order tensor
$
    \mathbf{U} = \Big[ 
    \mathbf{U}^{(1)}, ...,  \mathbf{U}^{(N_\mu)}
    \Big]
    \in
    \mathbb{R}^{N_\mu \times (N_t + 1) \times N_u},
$
which constitutes the training dataset in this work. See Figure \ref{fig:dataset_generation} for a visual description of this dataset.

\subsection{GPLaSDI}
GPLaSDI is a parameterized, data-driven reduced order modeling method which consists of three main features: 1) an autoencoder to learn a nonlinear spatial compression of data to a latent space; 2) dynamics identification within the latent space using Sparse Identification of Nonlinear Dynamics (SINDy) \cite{brunton2016discovering}; 3) Gaussian process interpolation of latent space dynamical systems over the parameter space to achieve parameterization. Here, we briefly describe GPLaSDI. More details may be found in \cite{bonneville2024gplasdi}.

The autoencoder is composed of two neural networks: an encoder $\phi_e$ parameterized by weights and biases $\theta_e$ and a decoder $\phi_d$ parameterized by weights and biases $\theta_d$. The encoder spatially compresses the full space data matrix $\mathbf{U} \in \mathbb{R}^{N_\mu \times (N_t + 1) \times N_u}$ to a latent space data matrix $\mathbf{Z} \in \mathbb{R}^{N_\mu \times (N_t + 1) \times N_z}$, 
where $N_z \ll N_u$. Conversely, the decoder spatially decompresses $\mathbf{Z}$ into a reconstructed full space representation $\widehat{\mathbf{U}} \in \mathbb{R}^{N_\mu \times (N_t + 1) \times N_u}$.
We define an autoencoder reconstruction loss as the mean squared error between $\mathbf{U}$ and $\widehat{\mathbf{U}}$,
\begin{equation} \label{eq:autoencoder_loss}
    \mathcal{L}_{AE}(\theta_e, \theta_d) = || \mathbf{U} - \widehat{\mathbf{U}} ||_2^2
    .
\end{equation}

$\mathbf{Z}$ may be interpreted as a collection of latent space trajectories $\mathbf{Z} = \Big[\mathbf{Z}^{(1)}, ...,  \mathbf{Z}^{(N_\mu)}\Big]$, where $\mathbf{Z}^{(i)} = \Big[  \mathbf{z}_0^{(i)}, ..., \mathbf{z}_{N_t}^{(i)} \Big]$ corresponds to parameter vector $\boldsymbol{\mu}^{(i)}$. In the SINDy interpretation, each trajectory is assumed to come from a dynamical system,
\begin{equation}
    \frac{d\mathbf{z}}{dt} 
    = 
    \mathbf{f}^{(i)}(\mathbf{z})
    ;
    \hspace{10pt}
    \mathbf{z}(0) = \mathbf{z}_0^{(i)}
    .
\end{equation}
The velocity $\mathbf{f}^{(i)}(\mathbf{z})$ is approximated by a linear combination of user-defined ``library" functions of $\mathbf{z}$,
\begin{equation}
    \mathbf{f}^{(i)}(\mathbf{z}) 
    = 
    \boldsymbol{\Phi(\mathbf{z})} \cdot \boldsymbol{\Xi}^{(i)T},
\end{equation}
where $\boldsymbol{\Phi(\mathbf{z})}$ is a matrix of candidate terms involving $\mathbf{z}$, and $\boldsymbol{\Xi}^{(i)}$ is an (as yet unknown) matrix of coefficients. In this work, we limit the library to constant and linear terms, thereby approximating the latent space dynamics with an affine ODE system. $\boldsymbol{\Xi}^{(i)}$ is determined by first computing latent space velocities $\dot{\mathbf{Z}}^{(i)}$ using a first-order finite difference scheme on $\mathbf{Z}^{(i)}$, then solving the linear regression problem
\begin{equation}
    \dot{\mathbf{Z}}^{(i)}
    =
    \boldsymbol{\Phi(\mathbf{Z}}^{(i)}) \cdot \boldsymbol{\Xi}^{(i)T}
    .
\end{equation}
Note that separate $\boldsymbol{\Xi}^{(i)}$ are learned for each $\boldsymbol{\mu}^{(i)}$. The collection of dynamics coefficient matrices is denoted $\boldsymbol{\Xi} = \Big[\boldsymbol{\Xi}^{(1)}, ..., \boldsymbol{\Xi}^{(N_\mu)} \Big]$. Associated with the dynamics identification is a loss term
\begin{equation} \label{eq:sindy_loss}
    \mathcal{L}_{SINDy} (\boldsymbol{\Xi})
    = \frac{1}{N_\mu} \sum_{i=1}^{N_\mu} 
    || 
    \dot{\mathbf{Z}}^{(i)} 
    - \boldsymbol{\Phi(\mathbf{z})} \cdot \boldsymbol{\Xi}^{(i)T} 
    ||_2^2
    .
\end{equation}
The model is trained by minimizing the following total loss function,
\begin{equation}
    \mathcal{L}(\theta_e, \theta_d, \boldsymbol{\Xi}) 
    = 
    \beta_1 \mathcal{L}_{AE}(\theta_e, \theta_d)
    +
    \beta_2 \mathcal{L}_{SINDy} (\boldsymbol{\Xi})
    +
    \beta_3 || \boldsymbol{\Xi} ||_2^2 
    ,
\end{equation}
where the third term is an important regularization term, and $\beta_1$, $\beta_2$, and $\beta_3$ are weighting hyperparameters.

The trained model can then be used to predict a time-dependent temperature field, given a new parameter vector $\boldsymbol{\mu}^{(*)}$ and a full space initial condition $\mathbf{u}_0^{(*)}$, which in our case is always uniform room temperature. First, $\mathbf{u}_0^{(*)}$ is compressed to the latent space using the encoder $\phi_e$, yielding a latent space initial condition $\mathbf{z}_0^{(*)}$. Then, a predicted latent space trajectory $\widetilde{\mathbf{Z}}^{(*)}$ is obtained by integrating in time, where the coefficient matrix $\boldsymbol{\Xi}^{(*)}$ is obtained by interpolating the $\boldsymbol{\Xi}^{(i)}$ from training samples over the parameter space. Gaussian process interpolation is used, which has the additional benefit of built-in uncertainty quantification. Finally, the decoder $\phi_d$ is used to generate the predicted full space trajectory $\widetilde{\mathbf{U}}^{(*)}$. Note that due to the uncertainty in $\boldsymbol{\Xi}^{(*)}$, multiple predicted trajectories produce slightly different results, and the mean and variance may be computed.

Another feature of GPLaSDI is a variance-based greedy sampling algorithm. During training, samples may be added to the training set every $N_{greedy}$ epochs. If there are currently $N_\mu$ samples in the training set, the additional $(N_\mu + 1)^{th}$ sample is chosen as that with the largest prediction uncertainty, as measured by the variance of the full order prediction.

\section{Results}
The following hyperparameters were found to yield satisfactory preliminary results. The encoder architecture follows a 25,990-1000-200-50-20-5 structure, comprising four hidden linear layers of size 1000, 200, 50, and 20, and a latent space of dimension 5. The decoder is symmetric with respect to the encoder. All layers except for the last in the encoder and decoder are connected with Softplus activation functions. The model is trained for $N_{epoch} = 1,000,000$ epochs using the Adam optimizer with a learning rate of $\alpha = 10^{-4}$. The training dataset begins with $N_\mu = 4$ samples corresponding to the corners of the parameter grid. Starting with corner points is an arbitrary but intuitive choice, since it ensures that predictions for all other points are interpolations rather than extrapolations.
A new sample is added to the training set every $N_{greedy}=200,000$ epochs. For the loss weights, we choose $\beta_1 = \beta_2 = 1$. The performance depends strongly on the value of $\beta_3$, and we compare results for $\beta_3$ of $10^{-3}$ and $10^1$. 

\begin{figure}[h!]
    \centering
    \includegraphics[width = 0.85\textwidth]{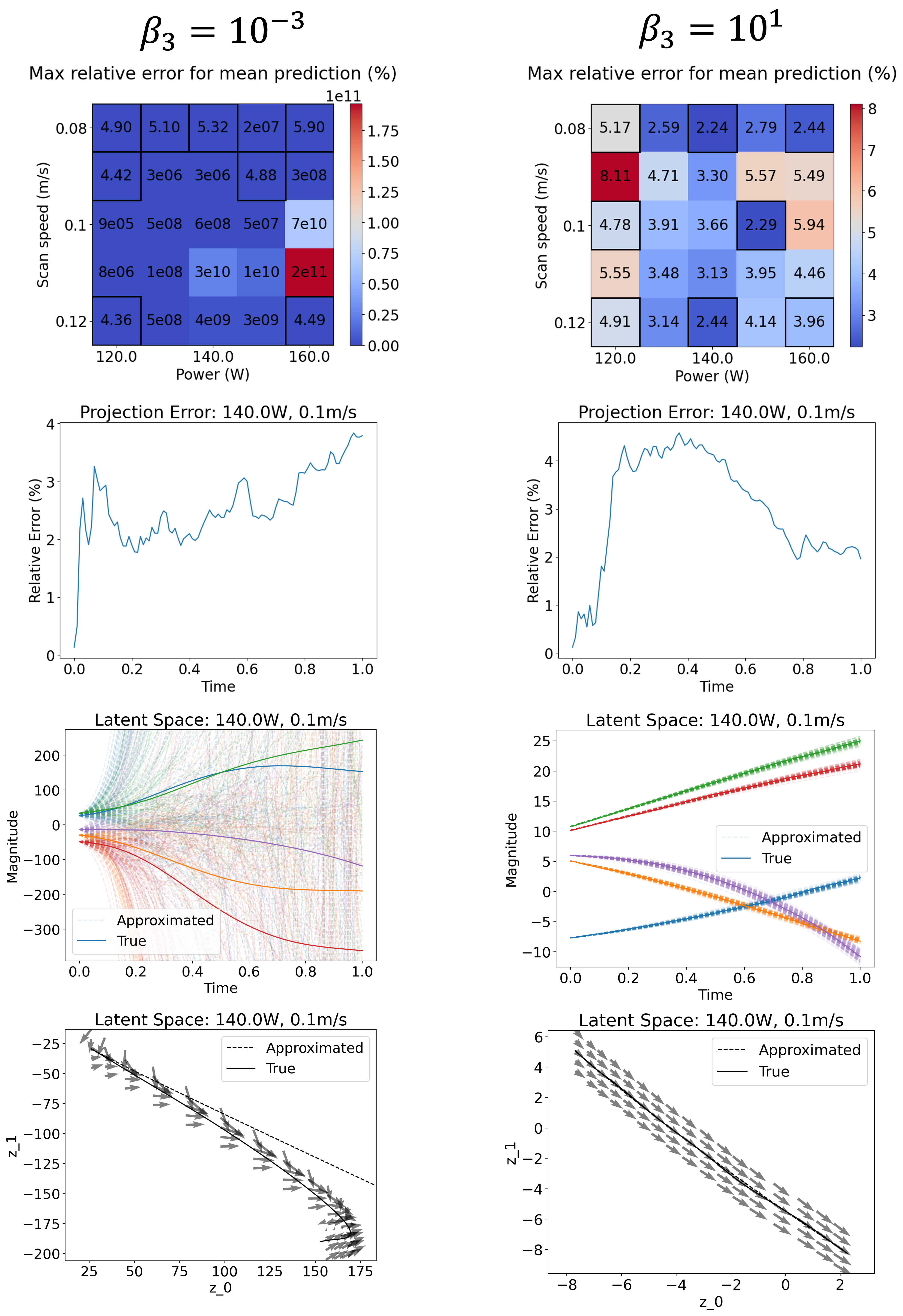}
    \caption{GPLaSDI results for $\beta_3 = 10^{-3}$ (left column) and $\beta_3 = 10^1$ (right column). First row: Maximum relative error for the mean GPLaSDI prediction over the parameter grid. Black squares denote training samples. Second row: Autoencoder projection error for $\boldsymbol{\mu}^{(*)} = (140 \mbox{W}, 0.1 \mbox{m/s})$. Third row: Latent space trajectories $\mathbf{Z}^{(*)}$ (true) and multiple $\widetilde{\mathbf{Z}}^{(*)}$ (approximated) for $\boldsymbol{\mu}^{(*)}$. Fourth row: Latent space trajectories (lines) and velocity (arrows) projected onto $z_0$-$z_1$ space for $\boldsymbol{\mu}^{(*)}$.
    }
    \label{fig:results} 
\end{figure}

To assess model performance, we use the maximum relative error, defined as
\begin{equation}
    e(\widetilde{\mathbf{U}}^{(*)}, \mathbf{U}^{(*)}) 
    = 
    \max_n 
    \Big( 
    \frac{||\widetilde{\mathbf{u}}_n^{(*)} - \mathbf{u}_n^{(*)}||_2}{|| \mathbf{u}_n^{(*)} ||_2}
    \Big).
\end{equation}
In Figure \ref{fig:results}, the left column shows the GPLaSDI results for $\beta_3 = 10^{-3}$, while the right column shows results for $\beta_3 = 10^{1}$. The top row shows the maximum relative error over the parameter grid, with training samples indicated by a black border. Clearly, $\beta_3$ is a very important hyperparameter, with extremely high errors for $\beta_3 = 10^{-3}$ and only $8.11\%$ test error in the worst case for $\beta_3 = 10^{1}$. In the second row, we compare the autoencoder projection errors vs. time, obtained by encoding and decoding the ground truth test data, for a typical test sample $\boldsymbol{\mu}^{(*)} = (140 \mbox{W}, 0.1 \mbox{m/s})$. For both values of $\beta_3$, projection errors are low, indicating that the autoencoder is not responsible for the vast difference in performance. Rather, the difference stems from disparate latent space dynamics behavior. In the third row, we show the latent space trajectories -- both the true $\mathbf{Z}^{(*)}$ and the (multiple) approximated $\widetilde{\mathbf{Z}}^{(*)}$ -- for the same typical test sample. For $\beta_3 = 10^{-3}$, the approximated trajectories are extremely unstable and blow up, while for $\beta_3 = 10^{1}$, the approximated trajectories match the true trajectory well. The stark difference in behavior arises from differences in the character of the dynamical system. In the fourth row, we visualize the same trajectories in a projected phase space, plotting the second latent space variable $z_1$ (orange lines) vs. the first latent space variable $z_0$ (blue lines). In addition, the arrows indicate the identified latent space velocity $\mathbf{f}^{(*)}(\mathbf{z}) = \boldsymbol{\Phi(\mathbf{z})} \cdot \boldsymbol{\Xi}^{(*)T}$ projected onto the $z_0$-$z_1$ plane. Arrows are drawn at points on and slightly offset from the true trajectory. As shown, the latent space velocity field for $\beta_3 = 10^{-3}$ is irregular and varies rapidly, promoting instability, while the velocity field for $\beta_3 = 10^{1}$ is very smooth, resulting in stable and accurate latent space trajectories. We note that increasing $\beta_3$ promotes smoothness because it tends to reduce the coefficients in $\boldsymbol{\Xi}^{(i)}$, which is closely related to $\nabla \mathbf{f}^{(i)}(\mathbf{z})$. With respect to computational cost, each GPLaSDI prediction takes on average $7.84 \mbox{ms}$ on an IBM Power9 CPU core with one NVIDIA V100 GPU, resulting in a roughly $1,000,000$ times speed-up compared to the physics-based simulation.

In GPLaSDI test cases \cite{bonneville2024gplasdi}, adequate performance was achieved with $\beta_3$ values of $10^{-5}$ or $10^{-6}$, much smaller than that required in our study. We suspect that the parameter grid density plays a important role, given the 21 $\times$ 21 grid in \cite{bonneville2024gplasdi} versus our 5 $\times$ 5 grid over a similar sized region of parameter space. From Figure \ref{fig:results}, we note that despite instability in identified dynamics, the poor-performing $\beta_3 = 10^{-3}$ model has relatively low error on training samples, leading us to attribute poor test sample predictions to the \textit{interaction} of unstable dynamics and interpolation over the parameter space. We propose that increasing grid density mitigates instability on test samples by reducing the interpolation challenge. Thus, our investigation of the role of $\beta_3$ may introduce a way to retain model performance with fewer training samples. Other factors, such as underlying physics and training data timestep size, may also influence the required $\beta_3$ size.

\section{Conclusion}
In this work, we applied GPLaSDI to time-dependent temperature data from single-track DED simulations for varying laser power and scan speed. With appropriate choice of hyperparameters, GPLaSDI is an effective reduced order model for this data, achieving $1,000,000$x speed-up compared to full-order, physics-based simulations with only about $8\%$ error in the worst case. In particular, the coefficient regularization weight $\beta_3$ has a large impact on model performance for the dataset considered in this work. For small $\beta_3$, learned latent space dynamics lead to unstable blow up of latent space trajectories on test samples. Increasing $\beta_3$ tends to smooth the latent space velocity field, preventing instabilities in the latent space trajectories and greatly improving prediction accuracy.

\begin{ack}
This work was performed under the auspices of the U.S. Department of Energy (DOE), by Lawrence Livermore National Laboratory (LLNL) under Contract No. DE-AC52–07NA27344. Funding: Laboratory directed research and development project 22-SI-007. IM release LLNL-CONF-854957.
\end{ack}


\bibliographystyle{unsrt} 
\bibliography{refs} 


\end{document}